\ifx\shlhetal\undefinedcontrolsequence\let\shlhetal\relax\fi

\documentstyle[11pt,amsfonts,amstex, amssymb]{amsart}
\newcount\skewfactor
\def\mathunderaccent#1#2 {\let\theaccent#1\skewfactor#2
\mathpalette\putaccentunder}
\def\putaccentunder#1#2{\oalign{$#1#2$\crcr\hidewidth
\vbox to.2ex{\hbox{$#1\skew\skewfactor\theaccent{}$}\vss}\hidewidth}}
\def\name{\mathunderaccent\tilde-3 }

\newcommand{\cD}{{\mathcal D}}
\newcommand{\Length}{{\rm Length}}
\newcommand{\Depth}{{\rm Depth}}
\newcommand{\bV}{{\bf V}}

\newcommand{\bP}{{\mathbb P}}
\newcommand{\bQ}{{\mathbb Q}}
\newcommand{\bR}{{\mathbb R}}
\newcommand{\At}{{\name{\cD}}}
\newcommand{\Dom}{{\rm Dom}}
\newcommand{\Levy}{{\rm Levy}}
\newcommand{\Reg}{{\rm Reg}}
\newcommand{\CON}{{\rm CON}}
\newcommand{\ZFC}{{\rm ZFC}}
\newcommand{\bd}{{\rm bd}}
\newcommand{\cf}{{\rm cf}}
\newcommand{\lesdot}{\mathrel{\mathord{<}\!\!\raise 
0.8 pt\hbox{$\scriptstyle\circ$}}} 

\newcommand{\QED}{\hfill\vrule width 6pt height 6pt depth 0pt\vspace{0.1in}} 
\newcommand{\Proof}{\noindent{\sc Proof} \hspace{0.2in}} 

\newtheorem{theorem}{Theorem}[section] 
 
\newtheorem{proposition}[theorem]{Proposition} 
 
\newtheorem{claim}{Claim}[theorem]

\theoremstyle{definition}
 
\newtheorem{definition}[theorem]{Definition}

\theoremstyle{remark}
\newtheorem{conclusion}[theorem]{Conclusion}

\newtheorem{remark}[theorem]{Remark}

\title{Length of Boolean algebras and ultraproducts}
\author{Menachem Magidor}
\address{Institute of Mathematics\\
The Hebrew University\\
Jerusalem 91904, Israel}
\email{menachem@@math.huji.ac.il}

\author{Saharon Shelah}
\address{Institute of Mathematics\\
The Hebrew University\\
Jerusalem 91904, Israel\\
and Rutgers University\\
Mathematics Department\\
New Brunswick, NJ 08854, USA}
\email{shelah@@math.huji.ac.il}

\date{\today} 

\thanks{The research was partially supported by ``Israeli Basic
Research Foundation''. Publication 433.}

\subjclass{Primary: 06Exx, 03E35; Secondary: 03E55} 

\begin{document}

\begin{abstract}
We prove the consistency with ZFC of ``the length of an ultraproduct of
Boolean algebras is smaller than the ultraproduct of the lengths''.
Similarly for some other cardinal invariants of Boolean algebras.
\end{abstract}

\maketitle

\section{Introduction} On the length of Boolean algebras (the cardinality of
linearly ordered subsets) see Monk \cite{M1}, \cite{M2} (and Definition
\ref{def1} below). In Shelah \cite[\S 1]{Sh:345} it is said that Koppelberg
and Shelah noted that by the {\L}o\'s theorem for an ultrafilter $\cD$ on
$\kappa$ and Boolean algebras $B_i$ ($i<\kappa$) we have 
\begin{enumerate}
\item[$(*)$]\quad  $|\prod\limits_{i<\kappa}\Length (B_i)/\cD|\leq
\Length (\prod\limits_{i<\kappa} B_i/\cD)$, and \\
$\mu_i<\Length (B_i) \quad\Rightarrow \quad |\prod\limits_{i<\kappa}
\mu_i/\cD |<\Length(\prod\limits_{i<\kappa} B_i / \cD)$.
\end{enumerate}
D.~Peterson noted that the indicated proof fails, but holds for regular
ultrafilters (see \cite{Pe97}).  Now the intention in \cite{Sh:345}
was for $\Length^+$, i.e.\  
\begin{enumerate}
\item[$(*)^+$]\quad $|\prod\limits_{i<\kappa}\Length^+(B_i)/\cD|\leq\Length^+(
\prod\limits_{i<\kappa} B_i/\cD)$,
\end{enumerate}
where $\Length^+(B)$ is the first cardinal not represented as the cardinality
of a linearly ordered subset of the Boolean Algebra (the only difference being
the case the supremum is not attained).

Here we prove that the statement $(*)$ may fail (see Theorem \ref{thm3} and
Proposition \ref{cl5}). The situation is similar for many cardinal invariants.

Of course, if $(*)$ fails then (using ultraproducts of $\langle ({\mathcal
H}(\chi), B_i): i<\kappa\rangle$ or see e.g.\ Ros{\l}anowski, Shelah
\cite[\S 1]{RoSh:534}) we have  
$\{i<\kappa:\Length^+(B)$ is a limit cardinal $\}
\in\cD$, and $\prod\limits_{i<\kappa}\Length^+(B_i)/\cD$ is $\lambda$--like
for some successor cardinal $\lambda$. Hence 
\[\{i<\kappa:\Length^+(B_i)\mbox{ is a regular cardinal}\}\in\cD\]
hence $\{i: \Length^+(B_i)\mbox{ is an inaccessible cardinal}\}\in D$,
so the example we produce is in some respect the only one possible.
(Note that our convention is that ``inaccessible'' means regular limit
($> \aleph_0$), not necessarily strong limit.) 

More results on cardinal invariants of ultraproducts of Boolean
algebras can be found in \cite{Sh:462}, \cite{Sh:479}, \cite{RoSh:534}
and \cite{Sh:620}, \cite{RoSh:651}. This paper is continued for other
cardinal invariants (in particular spread) in \cite{ShSi:677}.

We thank Otmar Spinas and Todd Eisworth for corrections and comments. 

\section{The main result}

\begin{definition}
\label{def1}
\begin{enumerate}
\item  For a Boolean algebra $B$, let its length, $\Length(B)$, be $\sup\{|X|:
X\subseteq B$, and $X$ is linearly ordered (in $B$) $\}$.
\item  For a Boolean algebra $B$, let its strict length, $\Length^+(B)$, be
$$
\sup\{|X|^+: X\subseteq B, \mbox{ and }X\mbox{ is linearly ordered (in
}B)\}.
$$
\end{enumerate}
\end{definition}

\begin{remark}
\begin{enumerate}
\item  In Definition \ref{def1}, $\Length^+(B)$ is (equivalently) the first
$\lambda$ such that for every linearly ordered $X\subseteq B$ we have $|X|<
\lambda$.
\item If \ $\Length^+(B)$ is a limit cardinal then
$\Length^+(B)=\Length(B)$; and if 
\ $\Length^+(B)$ is a successor cardinal then $\Length^+(B)=(\Length(B))^+$.
\end{enumerate}
\end{remark}

\begin{theorem}
\label{thm3}
Suppose $\bV$ satisfies {\rm GCH} above $\mu$ (for simplicity), $\kappa$ is
measurable,  $\kappa<\mu$, $\mu$ is $\lambda^+$-hypermeasurable (somewhat less
will suffice), $F$ is the function such that $F(\theta)=$ the first inaccessible
$>\theta$, and $\lambda=F(\mu)$ is well defined, and $\chi<\mu$,
$\chi>2^{2^\kappa}$. 

Then for some forcing notion $\bP$ not collapsing cardinals, except those in
the interval
$(\mu^+,\lambda)$ [so in $\bV^\bP$ we have $\mu^{++}=\lambda=F^\bV(\mu)$], and not
adding subsets to $\chi$, in $\bV^\bP$, we have:
\begin{enumerate}
\item[$(\alpha)$]  in $\bV^\bP$ the cardinal $\mu$ is a strong limit of
cofinality $\kappa$, 
\item[$(\beta)$] for some strictly increasing continuous sequence $\langle\mu_ i:
i<\kappa\rangle$ of (strong limit) singular cardinals $>\chi$ with limit
$\mu$, each $\lambda_i=: F^\bV(\mu_ i)$ is still inaccessible and for any
normal ultrafilter $\cD\in\bV$ on 
$\kappa$ we have: 
$$
\prod\limits_{i<\kappa}\lambda_ i/\cD\mbox{ has order type }
\mu^{++}=F^\bV(\mu).
$$
\end{enumerate}
\end{theorem}

\begin{definition}
\label{def3A}
A forcing notion $\bQ$ is directed $\mu$-complete if: for a directed 
quasi-order $I$ (so $(\forall s_0, s_1\in I)(\exists t\in I)(s_0\leq_I t \ \&\ s_1\leq
t)$) of cardinality $<\mu$, and $\bar p=\langle p_t: t\in I\rangle$ such
that $p_t\in \bQ$ and $s\leq_I t \Rightarrow p_s \leq_{\bQ} p_t$, there is $p\in
\bQ$ such that $t\in I \Rightarrow p_t\leq_{\bQ} p$.
\end{definition}

\Proof Without loss of generality for every directed $\mu$--complete
forcing notion $\bQ$ of cardinality 
at most $\lambda$ satisfying the $\lambda$-c.c., in $\bV^\bQ$ the cardinal $\mu$ is
still $\lambda$-hypermeasurable. [Why? If $\mu$ supercompact, use Laver
\cite{L}, if $\mu$ is just $\lambda$-hypermeasurable see more in Gitik Shelah
\cite{GiSh:344}.]  

Let $\bQ$ be the forcing notion adding $\lambda$ Cohen subsets to $\mu$, i.e.,\
$\{f: f$ a partial function from $\lambda$ to $\{0,1\}$, $|\Dom(f)|<\mu\}$.

In $\bV$, let $\bR=\Levy(\mu^+,<\lambda)=\{f:f$ a partial two place function
such that 
$[f(\alpha,i)$ defined\quad $\Rightarrow\quad 0<\alpha<\lambda\ \&\ i<\mu^+\
\&\ f(\alpha,i)<\alpha]$ and $|\Dom(f)|<\mu^+\}$ (so $\bR$ collapses all
cardinals in $(\mu^+,\lambda)$ and no others, so in $\bV^\bR$ the ordinal $\lambda$
becomes $\mu^{++}$). Clearly $\bR$ is $\mu^+$--complete and hence adds
no sequence of length $\le\mu$ of members of $\bV$. 

In $\bV^\bQ$, there is a sequence $\bar{\cD}=\langle\cD_i: i < \kappa\rangle$
of normal ultrafilters on $\mu$ as in \cite{Mg4} and, $\bar g=\langle g_{i,
j}: i<j \leq \kappa\rangle$, $g_{i, j} \in {}^\mu {\mathcal H}(\mu)$ witness
this (that is $\cD_i\in (\bV^\bQ)^\kappa / \cD_j$, in fact $\cD_i$ is equal to
$g_{i, j}/\cD_j$ in the Mostowski collapse of $(\bV^\bQ)^\kappa/\cD_j$). 
Let $\bar \At=\langle{\name{\cD}}_i:
i\leq\kappa\rangle$, $\bar{\name{g}}=\langle \name{g}_{i, j}: i<j\leq
\kappa\rangle$  be $\bQ$--names of such sequences. Note that a
$\bQ$--name $\name{A}$ of a subset of $\mu$ is an object of size $\leq\mu$,
i.e.,\ it consists 
of a $\mu$--sequence of $\mu$--sequences of members of $\bQ$, say $\langle\langle
p_{i, j}: j<\mu\rangle: i<\mu\rangle$, and function
$f:\mu\times \mu\longrightarrow\{0,1\}$ such that each $\{p_{i, j}: j<\mu\}$ is a
maximal antichain of $\bQ$ and $p_{i, j} \Vdash_\bQ ``i\in A \Leftrightarrow
f(i, j)=1$''. 
So the set of members of $\bQ$ and the
set of $\bQ$--names of subsets of $\mu$ are the same in $\bV$ and in $\bV^\bR$.  
So in $\bV^{\bR\times\bQ}$ the sequence  $\bar{\cD}$ still gives a sequence of
normal ultrafilters as required in \cite{Mg4} as witnessed by $\bar g=\langle
g_{i,j}: i<j\leq \kappa\rangle$. Also the Magidor forcing $\bP(\bar{\cD}, \bar
g)$ (from there) for changing the cofinality of $\mu$ to $\kappa$ (not
collapsing cardinals not adding  subsets to $\chi$, the last is just by fixing
the first element in the sequence) is the same in $\bV^\bQ$
and $\bV^{\bR\times\bQ}$ and has the same set of names of subsets of $\mu$. 
We now use the fact that $\bP(\bar{\cD}, \bar g)$ satisfies the
$\mu^+$-c.c. (see \cite{Mg4}). 
Let $\bP= (\bQ\times \bR)\ast \bP(\bar{\At}, \bar{\name{g}})$, so again every
$\bQ\ast \bP(\bar{\At}, \bar{\name{g}})$-name involves only $\mu$ decisions so
also $\bV^{\bQ\ast \bP(\bar{\At}, \bar{\name{g}})}$, $\bV^{(\bQ\times
\bR)\ast\bP(\bar{\At}, \bar{\name{g}})}$ have the same subsets of $\mu$.
So our only
problem is to check conclusion ($\beta $) of Theorem \ref{thm3}.

Let $\cD\in\bV$ be any normal ultrafilter on $\kappa$ (so this holds also in
$\bV^\bR$, $\bV^{\bR\times \bQ}$, $(\bV^{\bR\times \bQ})^{\bP(\bar{\At})}$). 

\begin{claim}
\label{cl}
In $\bV^{\bQ*
\bP(\bar{\At}, \bar{\name{g}})}$, the linear order
$\prod\limits_{i<\kappa}F(\mu_i)/\cD $ has 
true cofinality $F(\mu)=\lambda$.
\end{claim}

\par \noindent
{\em Proof of the Claim:} \qquad
Clearly $\Vdash_{\bQ}$ ``for $i<\kappa$ we have $\mu^\mu /
{\name{\cD}}_i$ is well ordered'' (as ${\name{\cD}}_i$ is
$\aleph_1$--complete) and 
$$
\begin{array}{ll}
\Vdash_{\bQ}& \mbox{ ``for }i<\kappa\mbox{ we have }\mu^\mu/
{\name{\cD}}_i\mbox{ has cardinality }2^\mu\\
\ & \mbox{ and even }\prod_{i< \mu} 2^{|i|} / {\name{\cD}}_i\mbox{ has
cardinality } 2^\mu\mbox{''}
\end{array}
$$
[Why? As $\mu=\mu^{<\mu}$ and $\cD_i$ is a uniform ultrafilter on $\mu$.
In details, let $h: {}^{\mu>}2 \rightarrow \mu$ be one-to-one, and for
each $\eta\in {}^\mu 2$ define $g_\eta\in {}^\mu\mu$ by $g_\eta(i)=
h(\eta\restriction i)$. Then
$$
\eta\neq \nu\in {}^\mu 2 \quad\Rightarrow\quad \{i<\mu:
g_\eta(i) = g_\nu(i)\}\mbox{ is a bounded subset of }\mu
$$ 
and hence its complement belongs to $\cD_i$ but $|\{g_\eta(i): \eta\in
{}^\mu 2\}| = 2^{|i|}$].

Consequently, for some $F^*\in {}^\mu\mu$ we have 
$$
\Vdash_\bQ ``\prod\limits_{\zeta<\mu} F^*(\zeta) /
{\name{\cD}}_i\mbox{ is isomorphic to }\lambda\mbox{''.}
$$

If we look at the proof in \cite{L} (or \cite{GiSh:344}) which we use
above, we see that w.l.o.g. $F^*$ is the $F$ above (and so does not
depend on $i$). 
So let $\name{f}_{i, \alpha}$ be $\bQ$-names such that  
$$
\begin{array}{ll}
\Vdash_{\bQ} &\mbox{ ``for }i<\kappa \mbox{ and }\alpha< \lambda,\quad
\name{f}_{i, \alpha} \in \prod_{\zeta<\mu} F(\zeta)\mbox{ and }\\
\ & \qquad f_{i, \alpha}/ {\name{\cD}}_i\mbox{ is the }\alpha\mbox{-th
function in }\prod\limits_{\zeta< \mu} F(\zeta)/ {\name{\cD}}_i\mbox{.''}
\end{array}
$$ 
In $\bV^{\bQ}$ let $\cD_\kappa= \bigcap\limits_{i< \kappa} \cD_i$ and
let $B_i\in \cD_i \setminus \bigcup\limits_{j< i} \cD_j$ be as in
\cite{Mg4} (you can also produce them straightforwardly), so $\langle
B_i: i< \kappa\rangle$ is a sequence of 
pairwise disjoint subsets of $\mu$. Define $f_\alpha\in {}^\mu\mu$
for $\alpha< \lambda$ as follows: 
$$
f_\alpha \restriction B_i = f_{i,
\alpha}\restriction B_i\ \mbox{ and }\ f_\alpha\restriction (\mu\setminus
\bigcup\limits_{i<\kappa} B_i)\mbox{ is constantly zero.}
$$ 
So $\langle
f_\alpha: \alpha< \mu\rangle$ is $<_{\cD_\kappa}$--increasing and
cofinal in $\prod\limits_{\zeta< \mu} F(\zeta)$. Let
${\name{\cD}}_\kappa$, ${\name{B}}_i$, ${\name{f}}_\alpha$ be
$\bQ$--names forced to be as above. Then $\Vdash_{\bQ}$ ``for $\alpha<
\beta< \lambda$ the set ${\name{A}}_{\alpha, \beta} = \{\zeta< \mu:
f_\alpha(\zeta) < f_\beta(\zeta)\}$ belongs to ${\name{\cD}}_\kappa$''.

Now in $\bV^{\bQ}$, one of the properties of Magidor forcing
$\bP(\bar{\cD}, \bar{g})$ is that 
$$
\begin{array}{ll}
\Vdash_{\bP(\bar{\cD}, \bar g)}& \mbox{ ``for every }A\in
\cD_\kappa=\bigcap\limits_{i<\kappa} {\cD}_i\\
\ & \qquad \mbox{ for every }
i<\kappa \mbox{ large enough we have }{\name{\mu}}_i\in A\mbox{''}
\end{array}
$$
(where $\langle {\name{\mu}}_i: i<\kappa\rangle$ is the increasing
continuous $\kappa$--sequence cofinal in $\mu$ which $\bP(\bar{\cD},
\bar{g})$ adds).

Since for every $p\in \bP(\bar{\cD}, \bar{g})$, for some $q\geq p$ we
have (recall from \cite{Mg4} that $F^q(i)$ is the set which $q$
``says'' $\name{\mu}_i$ belongs to (when $q$ does not forces a value
to $\name{\mu}_i$))
$$
[F^q(i)\subseteq A_{\alpha, \beta}\mbox{ for every }i<\kappa\mbox{ large
enough],}
$$ 
hence in $\bV^{\bQ*\bP(\bar{\name{\cD}}, \bar{\name{g}})}$,
for $\alpha< \beta< \lambda$ the set $\{i< \kappa: \neg
(f_\alpha(\mu_i)< f_\beta(\mu_i))\}$ is bounded, i.e., $\langle \langle
f_\alpha(\mu_i): i< \kappa\rangle: \alpha< \lambda\rangle$ is
$<_{J^{\bd}_\kappa}$--increasing in $\prod\limits_{i<\kappa}
F(\mu_i)$.

On the other hand, in $\bV^\bQ$, assume $p \Vdash_{\bP(\bar{\cD},
\bar{g})} \mbox{``}\name{f} \in \prod\limits_{i<\kappa} F(\name{\mu}_i)\mbox{''}$.
W.l.o.g. $p$ forces that $F^*(\name{\mu}_i)< \name{\mu}_{i+1}$, and so
by \cite{Mg4} (possibly increasing $p$), we have, for some function $h$,
$\Dom(h)=\kappa$, $h(i)\in [i]^{<\aleph_0}$, that above $p$, we know:
$\name{f}(i)$ depends on the value of $\name{\mu}_j$ for $j\in
\{i\}\cup h(i)$. So we can define a function $f^*\in {}^\mu\mu$: 
$$
\begin{array}{ll}
f^*(\zeta)=\sup\{\gamma: & \mbox{ for some }i<\kappa, \mbox{ possibly }
\name{\mu}_i\mbox{ is }\zeta\mbox{ and } \\
\ & \gamma \mbox{ is a possible
value (above }p)\mbox{ of }\name{f}(i)\}.
\end{array}
$$
So $f^*(\zeta)< F(\zeta)$, hence (in $\bV^\bQ$) for some $\alpha$,
$f^*<_{\cD_\kappa} f_\alpha$ and consequently $p\Vdash
``f^*<_{J^{\bd}_\kappa}f_\alpha$''. So $\langle  f_\alpha: \alpha<
\lambda\rangle$ is, in $\bV^{\bQ*\bP(\bar{\name{\cD}}, \bar{\name{g}})}$,
$<_{J^{\bd}_\kappa}$--increasing and cofinal in
$\prod\limits_{i<\kappa}F(\mu_i)$ which is more than enough for \ref{cl}.

Note that \ref{cl} holds in $\bV^{(\bR\times \bQ)*\bP(
\bar{\name{\cD}}, \bar{\name{g}})}$ too (remember that any sequence
$\langle\alpha_i<F(\mu_ i):i<\kappa\rangle$ in 
$\bV^{(\bQ\times\bR)*\bP(\bar{\name{\cD}}, \bar{\name{g}})}$ is
bounded by a function $f\in\big( 
\prod\limits_{\theta\in \Reg\cap \mu}F(\theta)\big)^\bV$ (see
\cite{Mg4}) and also is 
itself in $\bV^{\bQ*\bP(\bar{\At}, \bar{\name{g}})}$). 

But why, if $\alpha_i< F(\mu_i)$ and   
$\langle\alpha_ i:i<\kappa\rangle\in\bV^{(\bR\times \bQ)*
\bP(\bar{\At}, \bar{\name{g}})}$,  do we have that
$\prod\limits_{i<\kappa}\alpha_i/\cD$ has cardinality 
$\leq\mu^+$ (this means $<\lambda$)?

It suffices to prove this inequality in the universe
$\bV_1=\bV^{\bQ\ast \bP(\bar{\At}, \bar{\name{g}})}$. 

Now $\bV_1^\kappa/\cD$ is well founded, hence there is an isomorphism 
from $\bV_1^\kappa/\cD$ onto a transitive class which we now call
$M$ and let $j$ be the isomorphism (= the Mostowski collapse). As
$\mu$ in $\bV^{\bQ* \bP(\bar{\At}, \bar{\name{g}})}$ is strong 
limit $>\kappa$, clearly $\alpha<\mu\quad\Rightarrow\quad 
j(\langle\alpha: i<\kappa\rangle/\cD)<\mu$; and as $\cD$ is normal,
and $\langle\mu_i:i<\kappa\rangle$ is increasing continuous, we have
$j(\langle\mu_i: 
i<\kappa\rangle/\cD)=\mu$. As $\alpha_i<F(\mu_i)$ we have (by the {\L}o\'s theorem):
\[
\begin{array}{ll}
M\models & \mbox{``}j(\langle\alpha_i:i<\kappa\rangle/\cD)\mbox{ is an ordinal
smaller than }\\
\ &\qquad\qquad\mbox{ the first inaccessible }{>}\mu\mbox{''}.
\end{array}
\]
But the property ``not weakly inaccessible" is preserved by extending the universe
(from $M$ to $\bV_1$). So we finish. \QED$_{\ref{thm3}}$

\begin{remark}
\begin{enumerate}
\item  The proof has little to do with our particular $F$. Assume
$F:\mu\longrightarrow\mu$ and we add 
\[(*)\quad F(\mu)=\lambda\quad\mbox{ for }\lambda=
\big(j(F)\big)(\mu),\quad j
\mbox{ an appropriate elementary embedding.}\]
Then all the proof of \ref{thm3} works except possibly the last
sentence, which use some absoluteness of the definition of $F$. 
\item We can also vary $\bR$.
\item Let $\lambda_i= F(\mu_i)$.
In $\bV^\prime=\bV^{(\bQ\times \bR)*\bP(\bar{\At}, \bar{\name{g}})}$
we have $\mu=\mu^{<\kappa}$,
moreover $(\forall\alpha<\lambda)[\mid\alpha\mid^{<\kappa}<\lambda]$ and $\mu$
is strong limit. Hence if in $\bV^\prime$, $\bP$ is a forcing notion satisfying the
$\kappa$-c.c.\ of cardinality $<\mu$, and $\cD^*$ is an ultrafilter on
$\kappa$ extending $\cD$ then (in $(\bV^\prime)^\bP$) we have:\quad
the ultraproduct 
$\prod\limits_{i<\kappa}\lambda_i/\cD^*$ is $\lambda$--like.
\end{enumerate}
\end{remark}

\begin{proposition}
\label{cl5}
Suppose $\langle\lambda_i:i<\kappa\rangle$ is a sequence of (weakly) inaccessible
cardinals $>\kappa$, $\cD$ an ultrafilter on $\kappa$, and the linear order
$\prod\limits_{i<\kappa}(\lambda_i, <)/\cD$ is $\lambda$--like, $\lambda$ regular.
\begin{enumerate}
\item There are Boolean algebras $B_i$ (for $i<\kappa)$ such that:
\begin{enumerate}
\item[(a)] $\Length(B_i)=\Length^+(B_i)=\lambda_i$,
\item[(b)] $\Length^+(\prod\limits_{i<\kappa}B_i/\cD)=\lambda$,
\item[(c)] if $\lambda=\mu^+$ then $\Length(\prod\limits_{i<\kappa}B_i/\cD)=\mu$.
\end{enumerate}
\item  Assume $\lambda_i$ is regular $>\kappa$, 
$|B_i|=\lambda_i$, $B_i=\bigcup\limits_{\alpha<\lambda} B_{i,
\alpha}$, $B_{i,\alpha}$ increasing continuous in $\alpha$ (the $B_i$,
$B_{i, \alpha}$ are Boolean algebras of course) and we have 
$$
(\forall\alpha<\lambda_ i)(|B_{i,\alpha}|<\lambda_ i)
$$ 
and 
\begin{enumerate}
\item[$(*)_0$] if $i<\kappa$, $\delta<\lambda_ i$, $\cf(\delta)=\kappa^+$
then:\quad  $B_{i,\delta}\lesdot B_i$ ($\lesdot$ is complete subalgebra sign)
i.e.\
\begin{enumerate}
\item[$(*)^1_{B_i, B_{i,\delta}}$] if $x\in B_i\setminus B_{i,
\delta}$, $x\neq 0$ then for some $y^*$ 
we have:
\begin{enumerate}
\item[(a)] $y^*\in B_{i,\delta}$ and $y^*\neq 0$ 
\item[(b)] for any $y\in B_{i,\alpha}$ such that $B_i\models$`` $y\leq y^*\ \&\
y\neq 0$'':
\[B_ i\models\mbox{`` }y\cap x\neq 0\ \&\ x-y\neq 0\mbox{ ''}\]
\end{enumerate}
\end{enumerate}
\end{enumerate}
{\em Then} $\Length^+(\prod\limits_{i<\kappa}B_ i/\cD)\leq\lambda$.
\end{enumerate}
\end{proposition}

\Proof (1)\qquad  Let for $i<\kappa$ and  $\alpha<\lambda_ i$, $B_\alpha^i$ be
the Boolean algebra generated freely by $\{x^i_\zeta:\zeta<\alpha\}$ except   
\begin{enumerate}
\item[$\otimes$] $x_\zeta^{i,\alpha}\leq x_\xi^{i,\alpha}$ for $\zeta<\xi<
\alpha$.
\end{enumerate}
Let $B_i$ be the free product of $\{B^i_\alpha:\alpha<\lambda_ i\}$ so $B_i$
is freely generated by $\{x_\zeta^{i,\alpha}:\alpha<\lambda_i,\
\zeta<\alpha\}$ except for $\otimes$. 

Let $B_{i,\beta}$ be the subalgebra of $B_i$ generated by $\{x_\zeta^{i,\alpha}:
\alpha<\beta,\zeta<\alpha\}$.

Now clause (a) holds immediately, and the inequality $\geq$ in clause
(b) holds by the {\L}o\'s theorem, 
and the other inequality follows by  part (2) of the proposition. Lastly
clause (c) follows. 
\medskip

\noindent (2)\qquad  W.l.o.g. the set of members of $B_i$ is
$\lambda_i$, and the set of elements of $B_{i, \alpha}$ is an initial segment.

Let $S_i=\{\delta:\delta<\lambda_ i$, the set of members of $B_{i,\delta}$ is
$\delta$ and $\cf(\delta)=\kappa^+\}$

Let $(B,<^*,S)=\prod\limits_{i<\kappa}(B_i,<_ i,S_i)/\cD$, with $<_i$ the
order on the ordinals $<\lambda_i$ ($\leq$ is reserved for the order in the
Boolean algebra). So $(|B|,<^*)$ is $\lambda$--like (where $|B|$ is
the set of elements of $B$).

Let $\langle y_i:i<\lambda\rangle$ be an $<^*$ increasing sequence of members
of $B$. Let 
\[\begin{array}{ll}
S^{\prime} \stackrel{\rm def}{=}\{\delta<\lambda: &\cf(\delta)=\kappa^+
\mbox{ and }\{y_ i:i<\delta\}\mbox{ has in }(|B|,<^*)\\
\ &\mbox{a least upper bound which we call $y^\delta$ and it belongs to $S\}$.}
  \end{array}\]
Now clearly (read \cite[\S 1]{Sh:420} if you fail to see; or assume
$2^\kappa<\mu$ and \cite[2.3 p. 269]{Sh:111}):
\begin{enumerate}
\item[$\oplus$] $S'\subseteq \lambda$ is stationary.
\end{enumerate}
Note: $\oplus$ is enough, as if $X\subseteq B$ is linearly ordered by $<^B$,
let $y_i\in X$ for $i<\lambda$ be pairwise distinct; as $<^*$ is
$\lambda$--like w.l.o.g. $\langle y_i: i<\lambda\rangle$ is
$<^*$--increasing, and let $S^\prime $ be as above. 
For each $\delta\in S^\prime$ apply $(*)^1_{B,B\restriction\{y:y<^*
y_\delta\}}$ from the assumption to $y_\delta$, $y^\delta$ (holds by
{\L}o\'s theorem) and get $y^*_i$. Then apply 
Fodor lemma, and get a stationary subset $S^2$ of $S^\prime$ and an element $y^*$
such that for every $i\in S^2$ we have $y^*_i=y^*$. Now the set of $y_i$ for
$i\in S^2 $ is independent (check or see
\cite[4.1]{Sh:92}). \QED$_{\ref{cl5}}$

So putting Theorem \ref{thm3} and Proposition \ref{cl5} together

\begin{conclusion}
Assume $\CON(\ZFC + $ some $\mu$ is $\lambda^+_1$--hypermeasurable, for
some strongly inaccessible $\lambda_1>\mu)$. Then it is consistent that
for some $\kappa$, and sequence $\langle B_i: i<\kappa\rangle$ of
Boolean algebras and ultrafilter $\cD$ on $\kappa$ we have, for some
$\lambda$:
\begin{enumerate}
\item[(a)] $\prod\limits_{i<\kappa} \Length(B_i)/\cD = \lambda^+$
\item[(b)] $\Length(\prod\limits_{i<\kappa} B_i/\cD)=\lambda$.
\end{enumerate}
\end{conclusion}

\begin{remark}
\begin{enumerate}
\item  We can say more on when ultraproducts of free products of Boolean
algebras has not too large length. 
\item  We can use the disjoint sum of $\langle B_{i,\alpha}:\alpha<\lambda_i\rangle$
instead.
\item  In the proof of \ref{cl5} we actually have
$\Depth^+(B_i)=\Length^+(B_i)$ and 
$\Depth^+(B)=\Length^+(B)$, and so similarly without +;\quad where 

$\Depth^+(B)=\sup\{|X|^+: X\subseteq B$ is well ordered$\}$

$\Depth(B)=\sup\{|X|: X\subseteq B$ is well ordered$\}$.

\noindent So the parallel of \ref{cl5} holds for Depth instead $\Length$.
\item Recal $c(B)$ is the cellularity of a Boolean algebra $b$, i.e.,
$\sup\{|X|: X$ is a set of pairwise disjoint non zero elements$\}$. 
If in the proof of \ref{cl5} defining $B^i_\alpha$ we replace  
\begin{enumerate}
\item[$\otimes$] $x_\zeta^{i,\alpha}\leq x_\xi^{i,\alpha}$ 
\end{enumerate}
by 
\begin{enumerate}
\item[$\otimes^4$] $x_\zeta^{i,\alpha}\cap x_ \xi^{i,\alpha}=0$
\end{enumerate}
{\it then}  $c(B_i)=\lambda_i=c^+(B_i)$, $c(B)=\mu^+$, $\lambda=c^+(B)$.

(Same proof.)
\item  We can also get the parallel to \ref{cl5} for the independence
number.
Let $B^i$ be the Boolean algebra generated by $\{x^{i, \alpha}_\zeta:
\zeta< \alpha$ and $\alpha< \lambda_i\}$ freely except
\begin{enumerate}
\item[$(\otimes_5)$] $x^{i, \alpha}_\zeta\cap x^{i, \beta}_\xi = 0$ if
$\alpha< \beta< \kappa$, $\zeta< \alpha$, $\xi<\beta$.
\end{enumerate}
Let $I_i$ be the ideal of $B^i$ which $\{x^{i, \alpha}_\zeta: \zeta<
\alpha$, $\alpha< \lambda_i\}$ generates. Clearly it is a maximal
ideal. Let $B_{i, \alpha}$ be the ideal of $B^i$ generated by $\{x^{i,
\beta}_\zeta: \zeta< \alpha\mbox{ and }\beta<\alpha\}$. Again
w.l.o.g. the universe of $B_i$ is $\lambda_i$ and let 
$$
C_i=\{\delta< \lambda_i: \mbox{ for }x\in B_i\mbox{ we have: }x<
\delta\mbox{ iff } x\in B_{i, \delta} \vee -x\in B_{i, \alpha}\}.
$$ 
It is a club of $\lambda_i$.

The $B_{i,\beta}$ are not Boolean subalgebras of $B_i$, just {\em
Boolean subrings}; now $(*)_0$ in proposition \ref{cl5} is changed
somewhat.  We will have $P^{B_i}= I_i$ 
and $(B,<^*,P^*)=\prod\limits_{i<\kappa}(B_i,<^*_i,P^{B_i})/\cD$.

We know:  
\begin{enumerate}
\item[$(\alpha)$] $P^*$ is a maximal ideal of $B$ ( by {\L}o\'s Theorem)
\item[$(\beta)$] if $i<\kappa$, $\delta<\lambda_i$ is a limit ordinal
and $\delta\in C_i$ {\it then} for any $x\in P^{B_i}$ there are $x_0<\delta$,
$x_0\in P^{B_i}$ and $x_1\in P^{B_i}$ disjoint to all members of $P^{B_i}$
which are $<_i\delta$ and $x=x_0\cup x_1$. Similarly for $B$ add if
you like $Q_i= C_i\subseteq \lambda_i$.
\end{enumerate}
\end{enumerate}
\end{remark}

\bibliographystyle{literal-plain}
\bibliography{lista,listb,listx,listf}

\end{document}